\documentclass[12pt]{article}
\usepackage[latin1]{inputenc}
\usepackage{amssymb,amsmath,amsfonts}
\newtheorem{theorem}{Theorem}
\newtheorem{proposition}[theorem]{Proposition}
\newtheorem{lemma}[theorem]{Lemma}

\newtheorem{corollary}[theorem]{Corollary}

\newtheorem{remark}[theorem]{Remark}

\newcommand{\R}{\mathbb{R}}

\newcommand{\Sf}{\mathbb{S}}
\newcommand{\C}{\mathbb{C}}
\newcommand{\Hy}{\mathbb{H}}
\newcommand{\spa}{\mbox{span}}

\newcommand{\Ima}{\mbox{Im }}

\newcommand{\grad}{\mbox{grad\,}}

\newcommand{\po}{{\hspace*{-1ex}}{\bf .  }}

\newcommand{\nap}{\nabla^{\perp}}
\newcommand{\gnap}{\hat\nabla^{\perp}}

\newcommand{\nab}{\tilde\nabla}
\newcommand{\nabar}{\bar\nabla}
\newcommand{\la}{\lambda}
\def\<{{\langle}}
\def\>{{\rangle}}

\def\n{\nabla}

\def\a{\alpha}

\def\bea{\begin{eqnarray*} }
\def\eea{\end{eqnarray*} }
\def\be{\begin{equation} }
\def\ee{\end{equation} }

\def\nap{\nabla^\perp}
\def\proof{\noindent\emph{ Proof: }}
\def\qed{\ifhmode\unskip\nobreak\fi\ifmmode\ifinner
\else\hskip5 pt \fi\fi\hbox{\hskip5 pt \vrule width4 pt
height6 pt  depth1.5 pt \hskip 1pt }}

\begin{document}

\title{The dual superconformal surface}
\author{M. Dajczer and T. Vlachos}
\date{}
\maketitle

\begin{abstract} It is shown that a superconformal surface with arbitrary codimension
in flat Euclidean space has a (necessarily unique) dual superconformal surface if and only 
if the surface is S-Willmore, the latter a well-known necessary condition to allow a dual
as shown by Ma \cite{ma}. Duality means  that both surfaces envelope the same central 
sphere congruence and are  conformal with the induced metric. Our main result is that
the dual surface to a superconformal surface can easily be described in parametric 
form in terms of a parametrization of the latter. Moreover,  it is shown that the starting 
surface is conformally equivalent, up to stereographic projection in the nonflat case, to 
a minimal surface in a space form (hence, \mbox{S-Willmore}) if and only if  either the dual
degenerates to a point (flat case) or the two surfaces are conformally equivalent 
(nonflat case).
\end{abstract}

A  surface  $f\colon M^2\to\R^{n+2}$ in Euclidean space with codimension 
$n\geq 2$ is called \emph{superconformal}  if at any point  the ellipse 
of curvature is a nondegenerate circle.  Recall that the \emph{ellipse of 
curvature}  at $p\in M^2$ is the ellipse in the normal space $N_f M$ of $f$ 
at $p$ given by
$$
{\cal E}(p)=\{\a_f(X,X)\,:\, X\in T_pM\;\;\mbox{and}\;\;|X|=1\},
$$
where $\a_f$ denotes the second fundamental form of $f$ with values in 
the normal bundle; see \cite{gr} and references therein for several 
facts on this concept whose study started almost a century ago due to 
the work of Moore and Wilson \cite{mw1}, \cite{mw2}.  

Superconformality is invariant under conformal transformations since 
the property of ${\cal E}(p)$ being a circle is invariant under conformal changes 
of the metric  of the ambient space. Hence, the results in this paper belong to 
the realm of conformal (Moebius) geometry of surfaces and can also be stated in 
terms of surfaces in a space form.   

It was shown by Rouxel \cite{ro} that superconformal surfaces in codimension two  
always arise in pairs $f,\tilde f\colon M^2\to \R^4$ of dual surfaces that 
induce conformal metrics on $M^2$ and envelop a common central sphere congruence. 
Recall that the \emph{central sphere congruence} (or \emph{mean curvature sphere 
congruence}) of an Euclidean surface with any codimension is the family of 
two-dimensional spheres that are tangent to the surface and have the same mean 
curvature vector as the surface at the point of tangency. 
The concept of central sphere congruence (called the \emph{conformal Gauss map} 
in a different context by Bryant \cite{br}) is central in conformal geometry 
and was extensively studied since the turn of the last century, fundamentally 
due to the work of Thomsen \cite{th} and Blaschke \cite{bl}; see \cite{hj} for 
a detailed discussion of this subject.

Rouxel also discovered that the surface of centers of the  central sphere 
congruence is a minimal surface of $\R^4$.  If $f$ is free of minimal points,  
the \emph{surface of centers} is the locus of centers of the spheres in 
the congruence,  thus parametrically described by the map $g\colon M^2\to\R^{n+2}$ 
given by 
$$
g=f+\frac{1}{|H|^2}H
$$
where $H$ denotes the mean curvature vector field of $f$. 

In this paper, we consider  superconformal surfaces in Euclidean space in arbitrary
codimension.  To no surprise, the case of codimension two is rather special 
and this has much to do with the minimality of the surface of centers.  
In fact, this property and the classical Weierstrass representation 
of minimal surfaces allowed  Dajczer and Tojeiro \cite{dt} to provide a complete 
local parametric representation of all superconformal surfaces in $\R^4$.   
Moreover, they showed that the dual to a  superconformal surface in codimension 
two reduces to a point if and only if the surface is conformally equivalent, i.e., 
congruent by a conformal diffeomorphism of $\R^4$, to a holomorphic curve in $\C^2$.

By  a \emph{dual} to a surface  $f\colon M^2\to\R^{n+2}$ we mean an immersion 
$\tilde f\colon M^2\to\R^{n+2}$ that induces a  conformal metric and possess a 
common central sphere congruence, that is, at each point of $M^2$ the sphere in 
the their centrals sphere congruences is the same. In fact, for convenience we 
allow the dual to reduce to a single point.

For a locally conformally substantial superconformal surface in codimension higher 
than two that carries a dual superconformal surface, it turns out that the surface  
of centers is \emph{never} minimal. 
A surface being \emph{locally conformally substantial} means that 
the image under $f$ of any open subset of $M^2$ is not contained in a proper affine 
subspace or a sphere in the ambient space $\R^{n+2}$.
This and the fact that in higher codimension superconformality 
is not longer such a strong assumption, make unlikely the goal to obtain, 
a complete parametric classification as in \cite{dt}.  
Nevertheless, it seems natural to expect for  some class of superconformal 
surfaces  the existence of a dual surface similar to the case considered by 
Rouxel. In fact, this turns out to be the case for the superconformal surfaces 
that are  S-Willmore.

The concept of S-Willmore was introduced by Ejiri \cite{ej} as a special class 
of Willmore surfaces.  Ma \cite{ma} showed that being S-Willmore is the condition 
for a surface to have a dual that, in fact, is unique.  For a complex coordinate 
$z=y_1+iy_2$ associate to local isothermal coordinates superconformality means 
that the complex line bundle spanned by $\a_f(\partial_z,\partial_z)$ is isotropic 
and S-Willmore that it is holomorphic with respect to the normal connection.

It is well-known \cite{ej} that minimal surfaces in space forms are the basic 
examples of \mbox{S-Willmore} surfaces.  Hence, the ``trivial" examples of superconformal 
S-Willmore surfaces in Euclidean space are the ones conformally equivalent to 
minimal superconformal surfaces in Euclidean space and the images under stereographic 
projection of the same class of surfaces in the sphere or hyperbolic space. 
Euclidean minimal superconformal surfaces are called  \emph{1-isotropic} and admit a 
Weierstrass type representation  given in \cite{DG2} based on results in \cite{Ch}.
In the spherical case, this class of surfaces has been studied in different contexts, 
see \cite{BPW}, \cite{Mi} and \cite{Vl}. 

There are plenty of ``non-trivial" examples of superconformal S-Willmore surfaces in 
Euclidean space. For instance,  the image under stereographic projection  of any super 
Willmore surface in an even dimensional sphere is a superconformal S-Willmore surface. 
The class of super Willmore  surfaces was introduced and classified 
by Ejiri \cite{ej} in terms of isotropic holomorphic curves in complex 
projective spaces.
\vspace{1ex}

Note that in conformal geometry we may assume, at least locally, that the 
mean curvature of a surface never vanishes by composing with a conformal 
diffeomorphism.

\begin{theorem}\po\label{main} Let $f\colon M^2\to\R^{n+2},\,n\geq 3$, be a
regular locally conformally substantial superconformal surface. Then $f$ 
has a dual superconformal surface if and only if it is S-Willmore. 
Moreover, the dual surface can be parametrized as
$$
\tilde f=f+\frac{2}{|H|^2}(H)^\Lambda,
$$
where $\Lambda$ is the normal subbundle of rank $n-2$ of the  
surface of centers perpendicular to the plane subbundle of the first normal 
bundle $N_1^f$ of $f$  orthogonal to the mean curvature vector and
$(H)^\Lambda$ denotes taking the $\Lambda$-component. 
Furthermore, up to conformal equivalence, we have the following cases:
\begin{itemize}
\item[(i)] The dual reduces to a 
single point if and only if $f$ is a minimal surface. 
\item[(ii)]  The dual is obtained by composing $f$ with an  inversion and a 
reflection with respect to its center if and only if 
$f$ is the image under stereographic projection of a minimal
surface in the sphere $\Sf^{n+2}$.
\item[(iii)] The dual is obtained by composing $f$ with an inversion
if and only if $f$ is the image under stereographic projection of a minimal
surface in the hyperbolic space $\Hy^{n+2}$.
\end{itemize}
\end{theorem}

 The necessity of the surface being S-Willmore in the theorem is due to 
Ma \cite{ma} as already mentioned. A submanifold being \emph{regular} 
(or nicely curved) means that
the first normal spaces, i.e., the normal subspaces spanned by the second 
fundamental form, have constant dimension  and thus form a subbundle of 
the normal bundle. Notice that any isometric immersion is regular
along the connected components of an open dense subset of the manifold,
hence in local submanifold theory, as is the case of this paper, regularity 
is just a minor technical assumption. Finally, we mention that
part $(i)$ is known (see Remark on p.\ 339 of \cite{ej}) 
but we were not able to find a proof.  
\medskip

Any superconformal surface  in codimension two is S-Willmore, thus  there is no need
of such requirement in that case. The codimension three case is still quite special 
as shown by the following result.  

\begin{theorem}\po\label{main2}  Any superconformal Willmore surface 
$f\colon M^2\to\R^5$ is S-Willmore.
\end{theorem}

The paper concludes with a proof of the main result in \cite{dt} by means of the approach 
we developed here.

\section{Preliminaries}

In this section, we first recall some basic properties of the ellipse of curvature 
of a surface and then briefly discuss the notions of superconformal and S-Willmore 
surface.
\vspace{1,5ex}

Let $f\colon M^2\to\R^{n+2}$, $n\geq 2$, stand for an isometric immersion of a
two-dimensional Riemannian manifold into Euclidean space. Denote by 
$\a_f\colon TM\times TM\to N_f M$ its second fundamental form taking values 
in the normal bundle.  

Given an orthonormal basis $\{X_1,X_2\}$ of the tangent space $T_pM$ at $p\in M^2$,
denote $\a_{ij}=\a_f(X_i,X_j)$, $1\leq i,j\leq 2$.
Then, for any unit vector $v=\cos\theta X_1+\sin\theta X_2$ we have 
\be\label{eq:avv}
\a_f(v,v)=H+\cos2\theta\,\xi_1+\sin2\theta\,\xi_2,
\ee
where
$\xi_1=\frac{1}{2}(\a_{11}-\a_{22})$, $\xi_2=\a_{12}$
and $H=\frac{1}{2}(\a_{11}+\a_{22})$ is the mean curvature vector of $f$ 
at $p$. Thus, when $v$ goes once around the unit tangent circle, the vector 
$\a_f(v,v)$ goes twice around the ellipse of curvature ${\cal E}(p)$ of $f$ at $p$ 
centered at $H$. Clearly ${\cal E}(p)$ degenerates  into a line segment or a point 
if and only $\xi_1$ and $\xi_2$ are linearly dependent, that is, at points where 
the normal curvature tensor $R^\perp$ vanishes. It follows from (\ref{eq:avv}) 
that ${\cal E}(p)$ is a circle if and only if for some (and hence any) 
orthonormal basis of $T_pM$ it holds that
$$
\<\a_{12},\a_{11}-\a_{22}\>=0
\mbox\;\;\;\mbox{and}\;\;\;|\a_{11}-\a_{22}|=2|\a_{12}|.
$$

The complexified tangent bundle $TM\otimes \mathbb{C}$ is decomposed into
the eigenspaces of the complex structure $J$, denoted by $T^{\prime }M$ and 
$T^{\prime \prime }M$, corresponding to the eigenvalues $i$ and $-i$. The
complex structure of $M^2$ is determined by the orientation and the induced metric.
The second fundamental form   can
be complex linearly extended to $TM\otimes \mathbb{C}$ with values in the
complexified vector bundle $N_f M \otimes \mathbb{C}$ and then decomposed
into its $(p,q)$-components, $p+q=2,$ which are tensor products of $p$ many
1-forms vanishing on $T^{\prime \prime }M$ and $q$ many 1-forms vanishing 
on $T^{\prime }M$.

Taking local isothermal coordinates $\{y_1,y_2\}$ and $z=y_1+iy_2$, we have that
the surface $f$ is superconformal if and only if the 
$(2,0)$-part of the second fundamental form is  isotropic, or equivalently, if the 
complex line bundle $\a_f(\partial_z,\partial_z)$ is isotropic.
A surface  $f\colon M^2\to\R^{n+2}$ is called \emph{S-Willmore} \cite{ej}, \cite{ma} 
when the complex line bundle $\a_f(\partial_z,\partial_z)$ is parallel in the normal 
bundle, that is, if   
$$
\nap_{\partial_{\bar{z}}}\a_f(\partial_z,\partial_z)\;\;
\mbox{is parallel to}\;\;\a_f(\partial_z,\partial_z).
$$

It is well-known that any S-Willmore surface is always Willmore \cite{ej} but the 
converse is not true (cf.\ \cite{dw}) unless the substantial codimension is $n=2$.  
A surface being Willmore or S-Willmore is invariant under conformal diffeomorphisms
of Euclidean space. Recall that a surface $f\colon M^2\to\R^{n+2}$ is called 
\emph{Willmore} \cite{ej} if its mean curvature vector field $H$ satisfies the 
Willmore surface equation obtained as the Euler-Lagrange equation of the Willmore 
functional, namely, if
\be\label{willmore}
\Delta^\perp H-2|H|^2H+\Sigma_{i,j=1}^2\<H,\a_{ij}\>\a_{ij}=0
\ee
where $\Delta^\perp$ is the Laplacian in $N_f M$ and $X_1,X_2$ is an orthonormal 
frame.

Using the Codazzi equation, it follows that
$$
\nap_{\partial_z} H =\frac{2}{ \rho^2}
\nap_{\partial_{\bar{z}}}\a_f(\partial_z,\partial_z),
$$
where $ds^2=\rho^2|dz^2|$ is the induced metric. Thus, the surface is
S-Willmore if and only if
$\nap_{\partial_z}H$ is parallel to $\a_f(\partial_z,\partial_z)$ or,
equivalently, if
\be\label{parallel}
\nap_{V} H\;\;\mbox{is parallel to}\;\;\a_f(V,V)
\ee
for any $V\in T'M$.

\section{The proofs}

We proceed with the proofs of the results stated in the introduction.  We caution
that several arguments contain simple but long computations denominated straightforward
that may be only sketched.\vspace{1,5ex}

 In the sequel we denote by $f\colon M^2\to\R^{n+2}$, $n\geq 2$, a regular 
locally substantial superconformal surface.  The latter assumption 
is that the image under $f$ of any open subset of $M^2$ is not contained in a 
proper affine subspace of the ambient space. Recall that \emph{regular} means 
that the first normal spaces have constant dimension  and thus form a subbundle of 
the normal bundle. The \emph{first normal space} $N_1^f$ of $f$ at 
$p\in M^2$ is the normal subspace spanned by the second fundamental form, i.e.,
$$
N_1^f(p)=\spa\{\a_f(X,Y):X,Y\in T_pM\}.
$$

Under the above assumptions, it is easy to see that second fundamental form of 
the surface has the shape
$$
A_{\xi_1}=\begin{pmatrix} 
\lambda_1+\mu& 0\\ 
0& \lambda_1-\mu&
\!\!\end{pmatrix},\;\;
A_{\xi_2}=
\begin{pmatrix} 
\lambda_2& \mu\\ 
\mu& \lambda_2
\end{pmatrix}
\mbox\;\;\mbox{and}\;\;A_{\delta}=\lambda I
$$
with respect to orthonormal frames $\{X_1, X_2\}$ of the tangent  bundle 
and $\{\xi_1, \xi_2,\delta\}$ of the first normal subbundle $N_1^f$.
Thus the mean curvature vector field of $f$ is 
$$
H=\lambda_1\xi_1+\lambda_2\xi_2+\lambda\delta.
$$
Notice that we cannot have that $\mu=0$  on an open subset of $M^2$ since, 
otherwise,  $f$ would be totally umbilical along that set and this contradicts
being substantial. In particular, the special case $\dim N_1^f=2$  (in particular, 
if $n=2$) can only occur if $\lambda=0$.

From the Codazzi equations for $\xi_1,\xi_2$ and $\delta$ 
we obtain, respectively, that
$$
X_1(\lambda_1)-X_1(\mu)=-2\mu\Gamma_2-\mu\gamma_2+\lambda_2\gamma_1-\lambda\psi_{11},\;
X_2(\lambda_1)+X_2(\mu)=2\mu\Gamma_1-\mu\gamma_1+\lambda_2\gamma_2-\lambda\psi_{21},
$$
$$
X_1(\mu)-X_2(\lambda_2)=2\mu\Gamma_2+\gamma_2(\lambda_1+\mu)+\lambda\psi_{22},\;
X_2(\mu)-X_1(\lambda_2)=2\mu\Gamma_1+\gamma_1(\lambda_1-\mu)+\lambda\psi_{12},
$$
and
\be\label{10}
X_1(\lambda)=\psi_{11}(\lambda_1-\mu)+\lambda_2\psi_{12}-\mu\psi_{22},\;
X_2(\lambda)=\psi_{21}(\lambda_1+\mu)+\lambda_2\psi_{22}-\mu\psi_{12},
\ee
where we used the following notations:
$$
\Gamma_i=\<\nabla_{X_i}X_i,X_j\>,\; i\neq j,\;\; 
\gamma_i =\<\nap_{X_i}\xi_1,\xi_2\>\;\;\mbox{and}\;\;
\psi_{ij}=\<\nap_{X_i}\delta,\xi_j\>,\;\; i,j=1,2.
$$
The first four equations yield
\bea
X_1(\lambda_1)-X_2(\lambda_2)\!\!\!&=&\!\!\!\lambda_1\gamma_2
+\lambda_2\gamma_1+\lambda(\psi_{22}-\psi_{11}),\\
X_2(\lambda_1)+X_1(\lambda_2)\!\!\!&=&\!\!\!-\lambda_1\gamma_1
+\lambda_2\gamma_2-\lambda(\psi_{21}+\psi_{12}).
\eea
Setting 
\be\label{9}
X_1(\lambda_1)=\lambda_2\gamma_1-\lambda\psi_{11}-\mu a_1,\;\;
X_1(\lambda_2)=-\lambda_1\gamma_1-\lambda\psi_{12}-\mu a_2,
\ee
for some smooth functions $a_1,a_2$, we obtain that
\be\label{8}
X_2(\lambda_1)=\lambda_2\gamma_2-\lambda\psi_{21}+\mu a_2,\;\;
X_2(\lambda_2)=-\lambda_1\gamma_2-\lambda\psi_{22}-\mu a_1,
\ee
\be\label{12}
X_1(\mu)=\mu(2\Gamma_2+\gamma_2-a_1),\;\;
X_2(\mu)=\mu(2\Gamma_1-\gamma_1-a_2).
\ee

The Codazzi equation for any $\eta\in(N_1^f)^\perp$ is equivalent to
$$
\<\nap_{X_1}\eta,H\>=\mu(\<\nap_{X_1}\eta,\xi_1\>
+\<\nap_{X_2}\eta,\xi_2\>),
\;\;\<\nap_{X_2}\eta,H\>=\mu(\<\nap_{X_1}\eta,\xi_2\>
-\<\nap_{X_2}\eta,\xi_1\>).
$$
Let $\{\eta_\a\}_{1\leq\a\leq n-3}$
denote an orthonormal frame of $(N_1^f)^\perp$ and set
$$
\psi_{ij}^\a=\<\nap_{X_i}\eta_\a,\xi_j\>,\;\; 1\leq i,j\leq 2.  
$$
If $\dim N_1^f=2$ then $1\leq\a\leq n-2$. In the sequel, we work 
the case $\dim N_1^f=3$, but most of the computations hold  
if $\dim N_1^f=2$. For simplicity,  we denote
$$
\psi_1=\psi_{11}+\psi_{22},\;\;\psi_2=\psi_{21}-\psi_{12}\;\;\mbox{and}\;\;
\psi^\a_1=\psi^\a_{11}+\psi^\a_{22},\;\;\psi^\a_2=\psi^\a_{21}-\psi^\a_{12}.
$$
It follows using (\ref{10}) and (\ref{9}) to (\ref{12}) that
\be\label{13}
\nap_{X_1}H=-\mu(a_1\xi_1+a_2\xi_2+\psi_1\delta 
+\Sigma_\a\psi^\a_1\eta_\a),\;
\nap_{X_2}H=\mu(a_2\xi_1-a_1\xi_2+\psi_2\delta 
+\Sigma_\a\psi^\a_2\eta_\a).
\ee

The locus of the centers of the  central sphere congruence  given by 
the map $g\colon M^2\to\R^{n+2}$ defined as
$$
g=f+r^2H,\;\;\; \mbox{where}\;\;\;r=1/|H|,
$$
satisfies
\be\label{19}
g_*Z= f_*(I-r^2A_H)Z+r^2\nap_ZH+Z(r^2)H
\ee
where
$$
A_H
=\begin{pmatrix} 
|H|^2+\lambda_1\mu& \lambda_2\mu\\ 
\lambda_2\mu& |H|^2-\lambda_1\mu\\
\end{pmatrix}. 
$$
Using that
$$
A^2_{H} = 2|H|^2A_H-(|H|^4-\mu^2\theta)I
$$
where
$\theta=\lambda_1^2+\lambda_2^2=|H|^2-\lambda^2$,
it follows that
$$
\<g_{*}Z,g_{*}Y\>=r^4\mu^2\theta\<Z,Y\>+r^4\<\nap_Z H,\nap_YH\>.
$$
\be\label{15} \mbox{Thus}\;
f\mbox\;\mbox{and}\; g\;\mbox{are conformal}\;\iff |\nap_{X_1}H|=|\nap_{X_2}H|
\mbox\;\;\mbox{and}\;\;\<\nap_{X_1}H,\nap_{X_2}H\>=0.
\ee

\begin{proposition}\label{facts}\po
The following facts are equivalent:
\begin{itemize}
\item[(i)] The immersion $f$ is S-Willmore.
\item[(ii)]  The immersions $f$ and $g$ are conformal and $\nap H\subset N_1^f$.
\item[(iii)]  $\nap H\subset\Ima(\a_f-\<\;,\;\>H)$.
\item[(iv)]
$\psi_1=\psi_2=0\;\;\mbox{and}\;\;\psi^\a_1=\psi^\a_2=0,\;1\leq\a\leq n-3$.
\end{itemize}
\end{proposition}

\proof  On one hand,
$$
\a_f(X_1-iX_2,X_1-iX_2)=2\mu(\xi_1-i\xi_2).
$$
On the other hand, we have from (\ref{13}) that
$$
\frac{1}{\mu}\nap_{X_1-iX_2}H=-(a_1+ia_2)(\xi_1-i\xi_2)
-(\psi_1+i\psi_2)\delta-\Sigma_\a(\psi^\a_1+i\psi^\a_2)\eta_\a,
$$
and it follows from (\ref{parallel}) that $(i)$ and $(iv)$ are equivalent.  

 From (\ref{13}) we see that the right hand side  of (\ref{15}) 
is equivalent to
\be\label{16}
\psi_1^2+\Sigma_\a\,(\psi^\a_1)^2
=\psi_2^2+\Sigma_\a\,(\psi^\a_2)^2\;\;\mbox{and}\;\;
\psi_1\psi_2
+\Sigma_\a\,\psi^\a_1\psi^\a_2=0,
\ee
and the remaining of the argument follows easily from 
(\ref{13}) to (\ref{16}).\qed

\begin{corollary}\po\label{contained}  
If $f$ is S-Willmore then $(N_1^f)^\perp\subset N_gM$.
\end{corollary}
\proof  We have from (\ref{19}) that $\<g_*Z,\eta_\a\>=0,\;1\leq\a\leq n-3$, 
for any $Z\in TM$.\vspace{1.5ex}\qed 

We now prove the second result stated in the introduction.
\vspace{1.5ex}

\noindent \emph{Proof of Theorem \ref{main2}:} The Ricci equation 
$$
\<R^{\perp}(X_1,X_2)H,\xi_j\>=\<[A_H,A_{\xi_j}]X_1,X_2\>,\;j=1,2,
$$
together with (\ref{12}) and (\ref{13}) yield for $j=1$ that
\be\label{R1}
X_1(a_2)+X_2(a_1)-2a_1a_2+a_1\Gamma_1+a_2\Gamma_2+\psi_{11}\psi_2
+\psi_{21}\psi_1=2\mu\lambda_2
\ee
and for $j=2$ that
\be\label{R2}
-X_1(a_1)+X_2(a_2)+a^2_1-a_2^2-a_1\Gamma_2+a_2\Gamma_1+\psi_{12}\psi_2
+\psi_{22}\psi_1=-2\mu\lambda_1. 
\ee
On the other hand, 
\bea
\<\Delta^{\perp}H,\xi_j\>\!\!\!&=&\!\!\!X_1\<\nap_{X_1}H,\xi_j\>
+X_2\<\nap_{X_2}H,\xi_j\>
-\<\nap_{X_1}H,\nap_{X_1}\xi_j\>-\<\nap_{X_2}H,\nap_{X_2}\xi_j\>\\
\!\!\!&-&\!\!\!\Gamma_1\<\nap_{X_2}H,\xi_j\>-\Gamma_2\<\nap_{X_1}H,\xi_j\>.
\eea
Using (\ref{12}), (\ref{13}), (\ref{R1}) and (\ref{R2}) we easily obtain
$$
\frac{1}{\mu}\<\Delta^{\perp}H, \xi_1\>=-2\lambda_1\mu 
-\psi_1^2+\psi_2^2\;\;\;\mbox{and}\;\;\;
\frac{1}{\mu}\<\Delta^{\perp}H, \xi_2 \>= -2\lambda_2\mu 
+2\psi_1\psi_2.
$$
Also, 
$$
\Sigma_{i,j=1}^2\<\a_f(X_i,X_j),H\>\a_f(X_i,X_j)
=2|H|^2H+2\mu^2(\lambda_1\xi_1+\lambda_2\xi_2).
$$
Now, we have from (\ref{willmore}) that $f$ is Willmore if and only if 
$\psi_1=0=\psi_2$, and the result follows from Proposition \ref{facts}. 
\qed

\begin{proposition}\label{prop5}\po Let $f\colon M^2\to\R^{n+2}$ be a substantial 
superconformal S-Willmore surface with $\dim N_1^f=2$. If $n\geq 3$ then $f$ is minimal.
\end{proposition}

\proof The same proof given in Proposition \ref{facts} that
parts $(i)$ and $(iv)$  are equivalent still holds if $\dim N_1^f=2$. 
Thus $\psi^\a_1=0=\psi^\a_2$ for $1\leq\a\leq n-2$.
On the other hand, the Codazzi equation for $\eta_\a$ is
$$
\psi^\a_{11}A_{\xi_1}X_2+\psi^\a_{12}A_{\xi_2}X_2=
\psi^\a_{21}A_{\xi_1}X_1+\psi^\a_{22}A_{\xi_2}X_1.
$$
We obtain that
$$
\lambda_2\psi^\a_{11}-\lambda_1\psi^\a_{12}=0\;\;\mbox{and}\;\;
\lambda_1\psi^\a_{11}+\lambda_2\psi^\a_{12}=0. 
$$
But $\theta\neq 0$ would give $\psi^\a_{ij}=0$, which is not possible. Thus
$f$ is minimal.\qed

\begin{proposition}\label{prop6}\po Let $f\colon M^2\to\R^{n+2}$ be a superconformal  
S-Willmore surface with $\theta=0$ and $\dim N_1^f=3$. Then $f$ is minimal 
inside a sphere in $\R^{n+2}$.
\end{proposition}

\proof 
From (\ref{9}) and  (\ref{8}) we obtain $\psi_{ij}=0=a_1=a_2.$
Then (\ref{10}) implies that $|H|$ is constant.  Since $H \neq 0$, then (\ref{13})
and Proposition \ref{facts} show that the umbilical direction $H$ is parallel  in 
the normal connection.\qed  

\begin{remark}\po
{\em It follows from (\ref{19}) that the locus of the centers of the central 
sphere congruence of a non-minimal surface
is a point if and only if the surface is minimal in a sphere.}
\end{remark}

In the sequel, we also assume that $f$ is S-Willmore with $\theta\neq 0\neq\lambda$ 
everywhere. Set
$$
h_1=r^2(\la_2\xi_1-\la_1\xi_2),\;h_2=r^2H-\frac{1}{\la}\delta\;\;\mbox{and}\;\;
h_j=\eta_{j-2},\; 3\leq j\leq n-1.
$$
Thus $N_1^f$ is spanned by orthogonal vectors
$$
N_1^f=\spa\{h_1,h_2,H\}\;\;\mbox{where}\;\;|h_1|^2=r^4\theta\;\;\mbox{and}\;\;|h_2|^2
=\frac{r^2\theta}{\la^2}.
$$

\begin{lemma}\po\label{17} 
The following equations hold:
\be\label{o1}
h_{1_{*}}=g_*\circ J +\omega_{11}h_1+\omega_{12}h_2+\Sigma_\a\omega_{1\a}h_\a ,
\ee
\be\label{o2}
h_{2_{*}}=g_*  +\omega_{21}h_1+\omega_{22}h_2+\Sigma_\a\omega_{2\a}h_\a,
\ee
\be\label{o3}
h_{\a_*}=-\frac{1}{|h_1|^2}\omega_{1\a}h_1-\frac{1}{|h_2|^2}\omega_{2\a}h_2
+\Sigma_\beta \omega_{\a\beta} h_\beta,\;\mbox{where}
\ee
$$
\omega_{11}=-\frac{\lambda}{\theta r^2}d(\lambda r^2),\;\;
\omega_{12}=\frac{\lambda r^2}{\theta} (J \grad(\lambda/r^2))^*,\;\;
\omega_{21}=-\frac{1}{\lambda r^2\theta}(J\grad{\lambda})^*,
$$
$$
\omega_{22}=-\frac{1}{\lambda r^2\theta}d\lambda,\;\;
\omega_{1\a}=r^2(A_\a\omega_1-B_\a\omega_2),\;\;
\omega_{2\a}=-\frac{1}{\la}(C_\a\omega_1+D_\a\omega_2),
$$
$$
\omega_{\a\beta}=\<\nap h_\a,h_\beta\>
$$
where $\omega_i=X_i^*$, $i=1,2$, $Z^*$ denotes the 1-form dual to $Z\in TM$.
Also, 
$$
C_\a=\<\nap_{X_1}\delta,h_\a\>,\;D_\a=\<\nap_{X_2}\delta,h_\a\>,\;
A_\a=\la_1\psi^\a_{12}-\la_2\psi^\a_{11}\;\mbox{and}\;
B_\a=\la_1\psi^\a_{11}+\la_2\psi^\a_{12}.
$$
\end{lemma}

\proof A straightforward computation of the derivatives in the ambient 
space yields
\bea
\nabar_{X_1}(\lambda_2\xi_1-\lambda_1\xi_2)
\!\!\! &=& \!\!\!{\mu}f_*(-\lambda_2X_1+\lambda_1X_2)
-(\mu a_2+\psi_{12}\lambda)\xi_1
+(\mu a_1+\psi_{11}\lambda)\xi_2\\
\!\!\! &+& \!\!\!X_2(\lambda)\delta+\Sigma_\a A_\a h_\a,\\
\nabar_{X_2}(\lambda_2\xi_1-\lambda_1\xi_2)
\!\!\! &=& \!\!\!{\mu}f_*(\lambda_1X_1+\lambda_2X_2)
-(\mu a_1-\psi_{11}\lambda)\xi_1
-(\mu a_2-\psi_{12}\lambda)\xi_2\\
\!\!\! &-& \!\!\!X_1(\lambda)\delta-\Sigma_\a B_\a  h_\a.
\eea
Another straightforward computation using (\ref{13}), (\ref{19}) and that
\be\label{20}
X_1(1/r^2)=-2\mu(a_1\lambda_1+a_2\lambda_2),\;\;\;
X_2(1/r^2)=2\mu(a_2\lambda_1-a_1\lambda_2)
\ee
gives
\bea
h_{1_{*}}X_1\!\!\! &=& \!\!\!g_*X_2-\frac{\lambda}{\theta} X_1(\lambda r^2)h_1
-\frac{\lambda r^2}{\theta}  X_2(\lambda/r^2)h_2+r^2\Sigma_\a A_\a h_\a,\\
h_{1_{*}}X_2\!\!\! &=& \!\!\!-g_*X_1-\frac{\lambda}{\theta} X_2(\lambda r^2)h_1
+\frac{\lambda r^2 }{\theta} X_1(\lambda/r^2)h_2-r^2\Sigma_\a B_\a h_\a.
\eea
Similarly, we have
\bea
h_{2_{*}}X_1\!\!\! &=&\!\!\!g_*X_1+\frac{1}{\lambda
r^2\theta}(X_2(\lambda)h_1-X_1(\lambda)h_2)
-\frac{1}{\la}\Sigma_\a C_\a h_\a,\\
h_{2_{*}}X_2\!\!\!&=&\!\!\!g_*X_2-\frac{1}{\lambda
r^2\theta}(X_1(\lambda)h_1+X_2(\lambda)h_2)
-\frac{1}{\la}\Sigma_\a D_\a h_\a,
\eea
and (\ref{o1}) and (\ref{o2}) follow.  The third equation is just the 
Weingarten formula.
\qed\vspace{1,5ex}

We decompose $h_1$ and  $h_2$ into its tangent and normal components to $g$, 
namely,
\be\label{hdecomp}
h_1=g_*Y+\eta,\;\;  h_2=g_*Z+\xi.
\ee 

\begin{lemma}\po\label{ht} It holds that
\be\label{11}
Y=J\grad_g\varrho\;\;\mbox{and}\;\;Z=-\grad_g\varrho, 
\ee
where $\varrho=r^2/2$ and $J$ denotes a complex structure in TM. 
\end{lemma}

\proof Let $u$ be the conformal factor between the metrics  induced 
by $g$ and $f$ on $M^2$, that is,  $\<\,,\,\>_g=u\<\,,\,\>_f$. 
 From (\ref{13}), we have
\be\label{21}
\nap_{X_1}H=-\mu(a_1\xi_1+a_2\xi_2)\;\;\mbox{and}\;\;\nap_{X_2}H
=\mu(a_2\xi_1-a_1\xi_2).
\ee
We obtain using (\ref{19}), (\ref{20}) and (\ref{21}) that
\bea
g_*Y\!\!\!&=&\!\!\! 
\frac{1}{u}g_*(\< {h}_1,g_* X_1 \>X_1+\< {h}_1,g_*X_2\> X_2)\\
\!\!\!&=&\!\!\!\frac{\mu r^4}{u}g_*((a_2\lambda_1-a_1\lambda_2)X_1 
+(a_1\lambda_1+a_2\lambda_2)X_2)\\
\!\!\!&=&\!\!\!\frac{r^4}{2u}g_*(X_2(1/r^2)X_1-X_1(1/r^2)X_2)\\
\!\!\!&=&\!\!\!-\frac{r^4}{2u}g_*J\grad_f (1/r^2)\\
\!\!\!&=&\!\!\!\frac{1}{u}g_*J\grad_f\varrho.
\eea
The computation of $g_*Z$ is similar.\qed

\begin{lemma}\po\label{suubundle} The vector fields $\eta$ and $\xi$ 
are linearly independent everywhere.
\end{lemma}

\proof Assume that $c_1\xi+c_2\eta =0$ with $(c_1,c_2)\neq 0$ at $x\in M^2$. 
Then,
$$
c_1h_1+c_2h_2=g_* ( c_1Y+c_2 Z  ).
$$
Thus $g_*(c_1Z_1+c_2Z_2)$ is normal to $f$ at $x$ and (\ref{19}) implies that 
$$
c_1Z_1+c_2Z_2\in\ker(I-r^2 A_H).
$$
Since $\det (I-r^2 A_{H}) =-r^4\theta^2\mu^2$, we conclude that $c_1Z_1+c_2Z_2=0$. 
Hence  $c_1h_1+c_2h_2=0$, and this is a contradiction. \vspace{1,5ex}\qed

We now consider the orthogonal decomposition 
$$
N_gM=P\oplus\Lambda,
$$
where $P= \spa\{\eta,\xi\}$ and $\Lambda=(N_1^f)^\perp\oplus L$ with $\dim L=1$.
\vspace{1ex}

We observe that $Y$ (and hence $Z$) cannot vanish. 
Otherwise, from Lemma \ref{ht}  it follows that $|H|$ is constant. From (\ref{20}) 
we obtain $a_1=a_2=0$. Since $f$ is S-Willmore, working as in the proof of 
Theorem \ref{main2}, we see that the Ricci equation implies that (\ref{R1}) and 
(\ref{R2}) still hold. These immediately yield $\lambda_1=\lambda_2=0$, which 
is a contradiction.

\begin{lemma}\po\label{parametrization} The surface $f$ can be parametrized
in terms of $g$ as
$$
f=g-g_*\grad_g\varrho-\rho\,\xi+\Omega w,
$$
$$
\mbox{where}\;\;\rho=|\grad_g\varrho|^2/|\xi|^2,\;\; 
\Omega= \sqrt{2\varrho-\rho(\rho+1)|\xi|^2}\;\;
\mbox{and}\;\; w=-(H)^\Lambda/|(H)^\Lambda|\in L.
$$
\end{lemma}

\proof Using (\ref{19}) 
 we have
\bea
r^2(H)^{{g_*}(T M)} 
\!\!\!&=&\!\!\! \frac{r^2}{u}g_*(\<H,g_*X_1\>X_1+\<H,g_*X_2\> X_2)\\
\!\!\!&=&\!\!\!-\frac{r^4}{2u}g_*\grad_f(1/r^2)\\
\!\!\!&=&\!\!\! g_*\grad_g\varrho.
\eea
 From (\ref{19}), (\ref{hdecomp}) and (\ref{11}), we obtain that
$$
\<H,\eta \> =-\<H,g_*Y\>=0\mbox\;\;\mbox{and}\;\;
\<H,\xi\>=-\<H,g_*Z\>=\frac{1}{r^2}|\grad_g\varrho|^2.
$$
We also have $\<\eta,\xi\>=0$ from Lemma \ref{ht}, and the result follows.
\qed

\begin{lemma}\po\label{Hg} The mean curvature of $g$ satisfies
$$
(H_g)^{\spa \{\eta\}}=0\;\;\mbox{and}\;\; (H_g)^{\spa \{\xi\}}=-2|h_2|^{-2}\xi.
$$
\end{lemma}

\proof Using Lemma \ref{17}, we have 
$$
d\omega_{11}=-\lambda\theta^{-2}d(1/r^2)\wedge d\lambda
=\omega_{12}\wedge\omega_{21}.
$$
Computing $d^2h_1=0$ using (\ref{o1}) to (\ref{o3}) gives
\bea
0\!\!\!&=&\!\!\!d(g_*\circ J)+(g_*\circ J)\wedge\omega_{11}
+g_*\wedge \omega_{12}\\
\!\!\!&+&\!\!\!(d\omega_{12}-\omega_{11}\wedge\omega_{12}
-\omega_{12}\wedge\omega_{22}
+|h_2|^{-2}\Sigma_{\a}\omega_{1\a}\wedge\omega_{2\a})h_2\\
\!\!\!&+&\!\!\!\Sigma_{\a}(d\omega_{1\a}-\omega_{11}\wedge\omega_{1\a}
-\omega_{12}\wedge\omega_{2\a}
-\Sigma_{\beta}\omega_{1\a}\wedge\omega_{\a\beta})h_\a.
\eea
From
$$
\omega_{11}(X_1)X_1+\omega_{11}(X_2)X_2
=-\frac{\lambda}{\theta r^2}\grad(\lambda r^2),\;\;
\omega_{12}(X_2)X_1-\omega_{12}(X_1)X_2
=\frac{\lambda r^2}{\theta}\grad(\lambda/r^2)
$$
we obtain that
$$
(g_*\circ J)\wedge \omega_{11} + g_*\wedge \omega_{12}
=\frac{4\lambda^2}{\theta r^2} g_*Z\!*1,
$$
where $*1$ is the volume element.
Moreover, we have $d(g_*\circ J)=-2H_g\!*1$. We obtain 
\bea
2H_g \!*1\!\!\!\!&=&\!\!\!
\frac{4\lambda^2}{\theta r^2}\!*1\,g_*Z
+(d\omega_{12}-\omega_{11}\wedge\omega_{12}
-\omega_{12}\wedge\omega_{22}
+|h_2|^{-2}\Sigma_{\a}\omega_{1\a}\wedge\omega_{2\a})h_2\\
\!\!\!&+&\!\!\!\Sigma_{\a}(d\omega_{1\a}-\omega_{11}\wedge\omega_{1\a}
-\omega_{12}\wedge\omega_{2\a} -
\Sigma_{\beta}\omega_{1\a}\wedge\omega_{\a\beta})h_\a.
\eea
Hence $H_g$ is perpendicular to $\eta$ and
the $\xi$-component of  $H_g$  is $-(2\lambda^2/\theta r^2)\xi$.
\qed

\begin{proposition}\po If a superconformal S-Willmore surface $f\colon M^2\to\R^{n+2}$ 
is locally conformally substantial and free of minimal points,  then the locus 
of centers $g$  is  a minimal surface if and only $n=2$.
\end{proposition}

\proof If $n\geq 3$, it follows from Proposition \ref{prop5} and 
Proposition \ref{prop6} that $\dim N_1^f=3$ and that $\theta\neq 0$ on an open 
dense subset of $M^2$. Now, that $g$ is not minimal is a consequence of 
Lemma \ref{Hg}. The case $n=2$ follows from the result in \cite{ro}, 
i.e., our Theorem \ref{cod2}.
\vspace{1,5ex}\qed

We now prove the main result of this paper.\vspace{1,5ex}

\noindent \emph{Proof of Theorem \ref{main}:}
In view of Corollary \ref{contained} we may denote $n_\a=h_\a\in N_gM$ 
for $\a\geq 3$. Differentiating (\ref{hdecomp}) and using Lemma \ref{17}, gives
\be\label{g1}
\nabla_X Y-A_{\eta}X=JX+\omega_{11}(X)Y +\omega_{12}(X) Z,
\ee
\be\label{g2}
\nabla_X Z-A_{\xi}X=X+\omega_{21}(X) Y+\omega_{22}(X) Z,
\ee
\be\label{g3}
A_{n_\a}X=\omega_{1\a}(X)\frac{Y}{|h_1|^2}+\omega_{2\a}(X)\frac{Z}{|h_2|^2},
\ee
\be\label{g4}
\a_g(X,Y)+\gnap_X \eta=\omega_{11}(X)\eta 
+\omega_{12}(X)\xi +\Sigma_\a\omega_{1\a}(X)n_\a,
\ee
\be\label{g5}
\a_g(X,Z)+\gnap_X \xi=\omega_{21}(X) \eta
+\omega_{22}(X) \xi+\Sigma_\a\omega_{2\a}(X)n_\a,
\ee
\be\label{g6}
\gnap_Xn_\a=-\omega_{1\a}(X)\frac{ \eta}{|h_1|^2}
-\omega_{2\a}(X)\frac{\xi }{|h_2|^2}+\Sigma_\beta\omega_{\a\beta}(X) n_\beta,
\ee
where $\gnap$ denotes the induced connection in the normal bundle of $g$.

We claim that 
\begin{eqnarray}\label{several}
\<A_\xi X,Z\>-X(\rho)|\xi |^2-\frac{\rho}{2}X|\xi|^2
\!\!\!&=&\!\!\!\<A_\eta JX,Z\>-\rho\<\gnap_{JX}\xi,\eta\>\nonumber\\
\!\!\!&=&\!\!\!-\rho\<A_\xi X,Z\>-\<X,Z\>-\frac{1}{2}X|Z|^2.
\end{eqnarray}
 From the definition of the forms $\omega_{ij}$ we obtain
\be\label{omega}
\omega_{11}(X)+\omega_{12}(JX)=\frac{4\<Z,X\>}{|Z|^2+|\xi|^2}\;\;\;\mbox{and}
\;\;\;\omega_{21}(JX)=\omega_{22}(X).
\ee
Using (\ref{omega}) it follows from (\ref{g1}) and (\ref{g2}) that
$$
\<A_\xi Z,Z\>=\<\n_ZZ,Z\>-(\omega_{22}(Z)+1)|Z|^2,\;\;
\<A_\xi Y,Z \>=\<\n_YZ,Z\>-\omega_{22}(Y)|Z|^2,
$$
$$
\<A_\xi Y,Y \>=\<\n_Y Z,Y\>+(\omega_{22}(Z)-1)|Z|^2,\;\;
\<A_\eta Z,Z\>=-\<\n_ZZ,Y\>+\omega_{11}(Y)|Z|^2,
$$
$$
\<A_\eta Y,Z \>
=-\<\n_YZ,Y\>-(1+\omega_{11}(Z))|Z|^2+\frac{4|Z|^4}{|Z|^2+|\xi|^2},
$$
whereas from (\ref{g4}) and (\ref{g5}) that
$$
\<A_\xi Z,Z\>=-\<\gnap_Z\xi,\xi\>+\omega_{22}(Z)|\xi|^2,
$$
$$
\<A_\xi Y,Z\>=-\<\gnap_Y\xi,\xi\>+\omega_{22}(Y)|\xi|^2
=-\< \gnap_Z\eta,\xi\>-\omega_{11}(Y)|\xi|^2
$$
and
$$
\<A_\xi Y,Y\>=-\<\nap_Y\eta,\xi\>+
\omega_{11}(Z)|\xi|^2-\frac{4|Z|^2|\xi|^2}{|Z|^2+|\xi|^2}.
$$
Then, 
\be\label{s1}
\rho\<\gnap_Y\xi,\xi\>=-(\rho+1)\<A_\xi Y,Z\>+\<\n_Y Z,Z\>,
\ee
\be\label{s2}
\rho\<\gnap_Z\xi,\xi\>=-(\rho+1)\<A_\xi Z,Z\>+\<\n_Z Z,Z\>-|Z|^2,
\ee
\be\label{s3}
\rho\<\gnap_Z\xi,\eta\>=\rho\<A_\xi Y,Z\>+\<A_\eta Z,Z\>-\<\n_Z Y,Z\>,
\ee
\be\label{s4}
\rho \<\gnap_Y\xi,\eta\>= \rho\<A_\xi Y,Y\>+\<A_\eta Y,Z\>-\<\n_Y Y,Z\>+|Z|^2.
\ee
and
\be\label{24}
(2+\mathrm{tr} A_\xi)|Z|^2=\<\n_Z Z,Z\>-\<\n_Y Y,Z\>.
\ee
It is enough to argue for $X=Y$ and $X=Z$. We have using (\ref{11}) and
(\ref{s2}) that
\be\label{22}
\<A_\xi Y,Z\>-Y(\rho)|\xi|^2-\rho\<\gnap_Y\xi,\xi\>
=-\rho\<A_\xi Y,Z \>-\<\n_YZ,Z\>.
\ee
Similarly, from (\ref{11}) and (\ref{s2}) we also obtain
\be\label{23}
\<A_\xi Z,Z\>-Z(\rho)|\xi|^2-\rho\<\gnap_Z\xi,\xi\>
=-\rho\<A_\xi Z,Z\>-\<\n_ZZ,Z\>-|Z|^2,
\ee
and one equality in (\ref{several}) follows from (\ref{22}) and (\ref{23}).
Using (\ref{11}) and (\ref{s3}), we have
$$
\<A_\eta JY,Z\>-\rho\<\gnap_{JY}\xi,\eta\>
=-\rho\<A_\xi Y,Z\>-\<\n_Y Z,Z\>.
$$
Moreover, using (\ref{24}) and (\ref{s4}) we obtain
\bea
\<A_\eta JZ,Z\>\!\!\!&-&\!\!\!\rho\<\gnap_{JZ}\xi,\eta\>
=|Z|^2+\rho\<A_\xi Y,Y\>-\<\n_YY,Z\>\\
\!\!\!&=&\!\!\!-\rho\<A_\xi Z,Z \>+|Z|^2+\rho|Z|^2\mathrm{tr}A_\xi-\<\n_YY,Z\>\\
\!\!\!&=&\!\!\!-\rho\<A_\xi Z,Z\>
-\<\n_ZZ,Z\>+|Z|^2((\rho+1)\mathrm{tr}A_\xi+3).
\eea
Since Lemma \ref{Hg} gives $\mathrm{tr} A_\xi =-4/(\rho+1)$, we obtain
$$
\<A_\eta JZ,Z\>-\rho\< \gnap_{JZ} \xi,\eta\>
=-\rho\<A_\xi Z,Z\>-\<\n_ZZ,Z\>-|Z|^2,
$$
and this completes the proof of (\ref{several}).

Assume that $f$ is as in the statement. From Lemma \ref{parametrization} we have
$f=g+g_*Z-\rho\xi+\Omega w$. Then, define $f^-\colon M^2\to\R^{n+2}$ by
$$
f^-=g+g_*Z-\rho\xi+\Omega w_-
$$
where $w_-=-w$. We show that
$$
N_{f^-}M =\spa \{h_1, h_2, h, h_3, \dots h_{n-1} \}
$$
where $h=g_*\n r +(\rho/r)\xi-(\Omega/r)w_-$.
First  compute $f^-_*$ and use (\ref{g1}) to (\ref{g6}) to obtain
$$
\<f^-_*X,h_1\>=\<f^-_*X,h_2\> =\<f^-_*X,h_\a\>=0.
$$
To prove that also $h$ is normal to $f^-$ it is sufficient to see 
that $h$ is unitary and that $f^-=g-rh$. Observe also that (\ref{o3}) implies that
$N_1^{f^-}=\spa \{h_1, h_2, h\}$ and that the shape operators
of $f^-$ satisfy $A^-_{h_\a}=0$.

To prove that $f^-$ is superconformal we need to show that there
exist an orthogonal tangent basis $X_1, X_2=JX_1$ and functions $a,b$ such that
\be\label{c1}
(h_{1 _*} X_1)^{f^-_*(TM)}=af^-_*X_1+bf^-_*X_2,
\ee
\be\label{c2}
(h_{1 _*} X_2)^{f^-_*(TM)}=bf^-_*X_1-af^-_*X_2,
\ee
\be\label{c3}
(h_{2 _*} X_1)^{f^-_*(TM)}=-bf^-_*X_1+af^-_*X_2,
\ee
\be\label{c4}
(h_{2 _*} X_2)^{f^-_*(TM)}=af^-_*X_1+bf^-_*X_2,
\ee
\be\label{c5}
(h_{*} X_1)^{f^-_*(TM)}=\frac{1}{r}(-(1+b)f^-_*X_1+af^-_*X_2),
\ee
\be\label{c6}
(h_{*} X_2)^{f^-_*(TM)}=\frac{1}{r}(af^-_*X_1-(1-b)f^-_*X_2).
\ee
Using (\ref{o1})-(\ref{o3}), we see that
\bea
(h_{1 _*} X)^{f^-_*(TM)}=\frac{|\eta|^2}{|h_1|^2}g_*JX
-\frac{\<JX,Y\>}{|h_1|^2}\eta-\frac{|\eta|^2\<JX,Z\>}{|\xi|^2|h_1|^2}\xi
-\frac{\<JX,Z\>}{r^2}\Omega w_-,
\eea
\be\label{h2}
(h_{2 _*}
X)^{f^-_*(TM)}=\frac{|\eta|^2}{|h_1|^2}g_*X-\frac{\<X,Y\>}{|h_1|^2}\eta-
\frac{|\eta|^2\<X,Z\>}{|\xi|^2|h_1|^2}\xi-\frac{\<X,Z\>}{r^2}\Omega w_-
\ee
and
$$
(h_{*} X)^{f^-_*(TM)}=\frac{1}{r}(-f^-_*X+(h_{2 _*} X)^{f^-_*(TM)}).
$$
Thus
$$
(h_{1 _*} X_1)^{f^-_*(TM)}=(h_{2 _*} X_2)^{f^-_*(TM)} \;\;\mbox{and}\;\;
(h_{1 _*} X_2)^{f^-_*(TM)}+(h_{2 _*} X_1)^{f^-_*(TM)}=0.
$$
This means that (\ref{c1})-(\ref{c6}) are equivalent to (\ref{c3})
and (\ref{c4}), and we only have to choose the basis so
that (\ref{c3}) and (\ref{c4}) hold or, equivalently, that
\be\label{c3e}
f^-_*X_1=c(h_{2 _*} X_1)^{f^-_*(TM)}+d(h_{2 _*} X_2)^{f^-_*(TM)},\;
f^-_*X_2=d(h_{2 _*} X_1)^{f^-_*(TM)}-c(h_{2 _*} X_2)^{f^-_*(TM)}.
\ee
 From Lemma \ref{Hg} it follows that the self adjoint tensor field $L$ given by
$$
LX=X+\n_X Z+\rho A_\xi X-\Omega A_{w_-}X
$$
has zero trace. Let  $\{X_1, X_2=JX_1\}$ be  an orthonormal basis with
respect to the metric induced by $g$. Clearly, for suitable functions $c,d$, 
we have
\be\label{L}
LX_1=\frac{|\eta|^2}{|Z|^2+|\eta|^2}(cX_1+dX_2),\;\;
LX_2=\frac{|\eta|^2}{|Z|^2+|\eta|^2}(dX_1-cX_2),
\ee
which are actually the $g_*(TM)$-components of (\ref{c3e}).

Using (\ref{several}), (\ref{g4}) and (\ref{g5})
we see that the
$\xi$ and $\eta$ components of (\ref{c3e}) are equivalent~to
$$
\<LX_1, Z\>=\frac{|\eta|^2}{|Z|^2+|\eta|^2}\<cX_1+dX_2, Z\>,\;\;
\<LX_2, Z\>=\frac{|\eta|^2}{|Z|^2+|\eta|^2}\<dX_1-cX_2,Z \>,
$$
which are part of (\ref{L}).

The $w$-components of (\ref{c3e}) are equivalent to
$$
\<A_{w_-} X_1, Z\>+ \frac{X_1(\Omega)}{1+\rho}=
-\frac{\Omega}{r^2(1+\rho)}\<cX_1+dX_2,Z\>,
$$
$$
\<A_{w_-} X_2, Z\>+ \frac{X_2(\Omega)}{1+\rho}=
-\frac{ \Omega}{r^2(1+\rho)}\<dX_1-cX_2,Z \>.
$$
On account of (\ref{L}) and
$$
\frac{|\eta|^2}{|Z|^2+|\eta|^2}=\frac{\Omega^2}{r^2(1+\rho)}
$$
the above equations are equivalent to
$$
\<X+\n_X Z +\rho A_\xi X, Z\>+\frac{1}{1+\rho}\Omega X(\Omega)=0.
$$
We now argue that this holds. Indeed, this follows by differentiating
$|Z|^2(1+ \rho) + \Omega^2=r^2$
with respect to $X$ and using (\ref{g2}) and (\ref{g5}).

Finally we note that the $(N^f_1)^\perp$-components validity  of (\ref{c3e}) follows from 
$$
N_{f^-}M =\spa \{h_1, h_2, h, h_3, \dots h_{n-1} \},
$$ 
equation  (\ref{h2}) and Corollary \ref{contained}.

To complete the proof that $f^-$ is superconformal,
we need to show that the basis $\{X_1,X_2\}$
can be chosen to be orthonormal with respect to $g$.
An easy computation gives
$$
|(h_{2 _*}X_1)^{f^-_*(TM)}|^2=|(h_{2 _*}X_2)^{f^-_*(TM)}|^2
=\frac{|\eta|^2}{|Z|^2+|\eta|^2}
$$
and
$$
\<(h_{2 _*}X_1)^{f^-_*(TM)},(h_{2 _*}X_2)^{f^-_*(TM)}\>=0.
$$
Thus, in view of (\ref{c3e}) we have to show that
$$
|f^-_*X_1|^2=|f^-_*X_2|^2=\frac{|\eta|^2}{|Z|^2+|\eta|^2}(c^2+d^2)
\;\;\mbox{and}\;\;
\<f^-_*X_1,f^-_*X_2\>=0.
$$
Hence $f^-$ and $g$ are conformal and the desired orthonormal basis
with respect to the metric induced by $f^-$ is
$$
Y_j= \frac{1}{|\eta|}(|Z|^2+|\eta|^2)^{1/2}(c^2+d^2)^{1/2}X_j,\;\;j=1,2.
$$
In particular,  the mean curvature vector field of $f^-$ is given by
$H_-=(1/r)h$
and the locus of the centers of the corresponding central sphere congruence is
$$
f^-+\frac{1}{|H_-|^2}H_-=f^-+rh=g.
$$
To conclude the proof that $f^-$ is the dual to $f$ it remains to show that
$$
f_*(TM)\oplus\spa\{H\}= f_*^-(TM)\oplus\spa\{ H_-\}
$$
which follows from
$$
(f_*(TM)\oplus\spa\{H\})^\perp=\spa \{h_1, h_2, h_3, \dots h_{n-1}\}
=(f_*^-(TM)\oplus\spa\{H_-\})^\perp.
$$

Conversely, if $f$ allows a dual  superconformal surface then it is S-Willmore
by a result of Ma \cite{ma}. According to Lemma \ref{parametrization},
we have $f=g+g_*Z-\rho\xi+\Omega w$ and the dual to $f$ is given by
$$
\tilde f=g+g_*Z-\rho\xi-\Omega w.
$$
Then, the parametrization
$$
\tilde f=f+\frac{2}{|H|^2}(H)^\Lambda
$$
of the dual follows easily from $\Omega w=-r^2(H)^\Lambda$,
$\Lambda=(N_1^f)^\perp\oplus L$ and
$$
N_1^f=\spa\{h_1,h_2\}\oplus\spa\{H\}.
$$

Assume the dual reduces to a point $p_0$, i.e., 
$f=p_0 + 2\Omega w$.
On the other hand, from
$$
w^\perp= \<w,h_1\>\frac{h_1}{|h_1|^2}+
\<w,h_1\>\frac{h_2}{|h_2|^2}+\<w,H\>\frac{H}{|H|^2}
+\Sigma_\a\<w,\eta_\a\>\eta_\a
$$
and Lemma \ref{parametrization} we obtain that
$$
w^\perp=\<w,r^2(H)^{\Lambda}\> H=- \Omega H.
$$
Thus,
$$
(f-p_0)^\perp=\varphi H,\;\; \varphi=-2\Omega^2.
$$
Moreover, we have
$$
|f-p_0|^2+2\varphi=0.
$$
Consider the inversion $\mathcal{I}$  with respect to a sphere with 
radius $R=1$ centered at $p_0$ and the immersion 
$\tilde{f}=\mathcal{I}\circ f$.  Then, there is  a vector bundle isometry
$\mathcal{P}$ between the normal bundles $N_fM$ and $N_{\tilde{f}}M$
(see \cite{dt0}) given by
$$
\mathcal{P}\mu=\mu-2\frac{\<f-p_0,\mu\>}{|f-p_0|^2}(f-p_0)
$$
such that shape operators of $f$ and $\tilde{f}$ are related by
$$
\tilde{A}_{\mathcal{P}\mu}=|f-p_0|^2A_\mu+2\<f-p_0,\mu\>I .
$$
We can easily find that the mean curvature vector of $\tilde{f}$ is given by
$$
H_{\tilde{f}}=\mathcal{P}(|f-p_0|^2 H+2(f-p_0)^\perp).
$$
Using  $(f-p_0)^\perp=\varphi H$ and 
$|f-p_0|^2+2\varphi=0,$ we deduce that $\tilde{f}$ is  minimal in
$\R^{n+2}$.\vspace{1,5ex}

Conversely, assume that  the  surface is a composition of a $1$-isotropic
surface with  an
inversion with respect to the sphere in $\R^{n+2}$ with radius $R=1$ centered
at $p_0$. Then,
$$
(f-p_0)^\perp=\varphi H\;\;\mbox{and}\;\;\varphi=-\frac{1}{2}\< f-p_0,f-p_0\>.
$$
 From this we obtain that
$$
\a(X,\grad \varphi)=X(\varphi)H+\varphi \nap_X H.
$$
Using (\ref{19}) we see that $f-p_0$ is perpendicular to the surface $g$. Since
it is also perpendicular to $h_1, h_2$, it follows
that it is  perpendicular to the plane bundle $P$. Furthermore, it is
perpendicular to $(N_1^f)^\perp$.  Thus, we conclude that
$f-p_0=\sigma w$.  From this we obtain 
$$
\frac{1}{\sigma}(f-p_0)^\perp= \<w,h_1\>\frac{h_1}{|h_1|^2}+
\<w,h_1\>\frac{h_2}{|h_2|^2}+\<w,H\>\frac{H}{|H|^2}
+\Sigma_\a\<w,\eta_\a\>\eta_\a
$$
or, using Lemma \ref{parametrization}, that
$$
(f-p_0)^\perp=\sigma\<w, (r^2 H)^{\Lambda}\> H=-\sigma\ \Omega H.
$$
Hence $\varphi=-\sigma\Omega$. In view of $\sigma^2=-2\varphi$, we have that
$\Omega=|f-p_0|/2$ and $\sigma=|f-p_0|$. Thus
$f-p_0=2\Omega w$, which shows that the  dual to $f$ reduces to the point $p_0$.
Finally, the cases (ii) and (iii) follow directly from Proposition \ref{dual}
given next.\vspace{1,5ex}\qed

In the following result we do not assume that the surface is superconformal.

\begin{proposition}\po\label{dual}
 Let $f\colon M^2\to\R^{n+2}$ be a surface with dual surface $\tilde f=T \circ f$
where $T$ denotes a conformal diffeomorphism of $\R^{n+2}$. 
\begin{itemize}
\item[(i)] Then $T$ is composition of an inversion with a reflection with respect 
to the center of the inversion if only if $f$ is either a composition of a minimal surface 
in the sphere \mbox{$\Sf^{n+2}\subset\R^{n+3}$}  with a  stereographic projection 
onto $\R^{n+2}$ or a minimal surface in a sphere $\Sf^{n+1}\subset\R^{n+2}$. 
\item[(ii)] Then $T$ is an inversion if and only if $f$ is a composition of a minimal 
surface in hyperbolic space $\Hy^{n+2}\subset\mathbb{L}^{n+3}$ with a stereographic 
projection onto $\R^{n+2}$. 
\end{itemize}
\end{proposition}

\proof A surface $\tilde f\colon M^2\to\R^{n+2}$ is the dual to a surface 
$f\colon M^2\to\R^{n+2}$ if 
\be\label{surfcen}
f+r^2H=\tilde f+r^2\tilde H,
\ee
\be\label{eqspe}
f_*(TM)\oplus\spa\{H\}=\tilde f_*(TM)\oplus\spa\{\tilde H\}
\ee
and $|H|=|\tilde  H|=1/r$.

To prove $(i)$ first assume that $\tilde{f}=-\mathcal{I}\circ f$, where 
$\mathcal{I}$  is the inversion with respect to a sphere  with radius $1$ 
centered at the origin.  The mean curvature vector field $\tilde H$ is given 
by $\tilde H=-\hat H$, where $\hat H $ is the mean curvature of the 
surface $\hat f=\mathcal{I}\circ f=f/|f|^2$. From the results in \cite{dt0} 
there is a vector bundle isometry $\mathcal{P}$ between the normal bundles 
$N_fM$ and $N_{\hat{f}}M$ given by
$$
\mathcal{P}\mu=\mu-\frac{2}{|f|^2}\<f,\mu\>f
$$
such that shape operators of $f$ and $\hat{f}$ are related by
$$
\hat{A}_{\mathcal{P}\mu}=|f|^2A_\mu+2\<f,\mu\>I 
$$
and the mean curvature vectors by
$$
\hat{H}= \mathcal{P}(|f|^2H+2f^\perp).
$$
Using (\ref{surfcen}) we deduce that 
$$
\big(1-2r^2\<f,H\>+\frac{1}{|f|^2}(1-4r^2|f^\perp|^2)\big)f
=-r^2((1+|f|^2)H+2f^\perp).
$$
Thus, the left hand side must vanish unless $f\in N_fM$. In the latter case
$f$ is a minimal surface in a sphere. If not we have 
$$
2f^\perp=-(1+|f|^2)H.
$$
Let $e$  be a unit vector in 
$\R^{n+3}=\R^{n+2}\oplus\R e$ and let $\mathcal{T}$ be the inversion 
$$
\mathcal{T}(p)=q_0+\frac{1}{|p-q_0|^2}( p-q_0)
$$
with respect to the sphere $\mathbb{S}^{n+2}$ with 
radius $1$ centered at $q_0=(0,1)=e$. If $\bar{f}=\mathcal{T}\circ f$,  there is,  
as before, a vector bundle isometry $\mathcal{\bar P}$ between $N_fM$ and 
$N_{\bar{f}}M$ given by
$$
\mathcal{\bar P}\mu=\mu-\frac{2}{|f-q_0|^2}\<f-q_0,\mu\>(f-q_0)
$$
such that shape operators of $f$ and $\bar{f}$ are related by
$$
\bar{A}_{\mathcal{\bar P}\mu}=|f-q_0|^2 A_\mu+2\<f-q_0,\mu\>I
$$
and the mean curvature vector field  of $\bar{f}$  is given by
\bea
\bar H\!\!\!&=&\!\!\!\mathcal{\bar P}(|f-q_0|^2 H+2(f-q_0)^\perp)\\
\!\!\!&=&\!\!\!\mathcal{\bar P}(|f-q_0|^2H-(1+|f|^2)H-2e)\\
\!\!\!&=&\!\!\!-2 \mathcal{\bar P}e\\
\!\!\!&=&\!\!\!-4(\bar f-(1/2)e).
\eea
Thus $\bar{f}$ is minimal in the sphere $\Sf_{1/2}^{n+2}(e/2)$ with radius $1/2$ 
centered  at $e/2$. 

Conversely, let $f$ be a composition of a minimal surface in the sphere  
$\Sf_{1/2}^{n+2}(e/2)$ in  $\R^{n+3}=\R^{n+2}\oplus\R e$ 
with a  stereographic projection onto $\R^{n+2}$. Then, 
$$
f^\perp= -\frac{1}{2}|f-q_0|^2H=-\frac{1}{2}(1+|f|^2)H
$$
where $q_0=e$. Consider the surface $\tilde{f}=-\mathcal{I}\circ f$, where  
$\mathcal{I}$ is the inversion with respect to a sphere with radius $1$ 
centered at the origin. 
We claim that $\tilde f$ is the dual to $f$. Consider the surface 
$\hat f= f/|f|^2$ with mean curvature $\hat H$ and $\mathcal P$ the corresponding 
bundle  isometry between the normal bundles of $f$ and $\hat f$.
The mean curvature of $\tilde f$ is given by
$$
\tilde H=-\hat H=-\mathcal{P}(|f|^2 H+2f^\perp)=\mathcal{P}(H).
$$
We have that $|\tilde H|=|H|$ and 
\bea
\tilde{f}+r^2 \tilde H \!\!\!&=&\!\!\!-\frac{f}{|f|^2}+r^2\mathcal{P}(H)\\
\!\!\!&=&\!\!\!-\frac{f}{|f|^2}+r^2 H-\frac{2r^2}{|f|^2}\<f^\perp,H\>f\\
\!\!\!&=&\!\!\!-\frac{f}{|f|^2}+r^2 H+\frac{1+|f|^2}{|f|^2}f\\
\!\!\!&=&\!\!\!f+r^2 H.
\eea
It remains to show that (\ref{eqspe}) holds  or, equivalently, that
$$
f_*(TM)\oplus\spa\{H\}=\mathcal{P}(f_*(TM)\oplus\spa\{H\}),   
$$
since 
$$
\tilde f_*=-\frac{1}{|f|^2} \mathcal{P}\circ f_*.
$$
Now decompose $f$ in its tangent and normal components as $f=f_*V+f^\perp$.
Then,
\bea
\mathcal{P} H \!\!\!&=&\!\!\! H-\frac{2}{|f|^2}\<f^\perp,H\>f\\
\!\!\!&=&\!\!\!-\frac{2}{|f|^2}\<f,H\>f_*V
+\Big(1+\frac{1}{|f|^2}\<f,H\>(1+|f|^2)\Big)H\in f_*(TM)\oplus\spa\{H\}.
\eea
Moreover, 
\bea
\mathcal{P} f_*X \!\!\!&=&\!\!\! f_*X- \frac{2}{|f|^2}\<f,f_*X \>f\\
\!\!\!&=&\!\!\!f_*(X-2\<X,V\>V)-\frac{2}{|f|^2}\<X,V\>f^\perp\\
\!\!\!&=&\!\!\!f_*(X-2\<X,V\>V)
+\frac{1}{|f|^2}\<X,V\>(1+|f|^2))H\in f_*(TM)\oplus\spa\{H\}.
\eea
Thus $\tilde f$ is the dual to $f$.
\vspace{1,5ex}

To prove $(ii)$, we proceed as in the previous case assuming that  
$\tilde{f}=\mathcal{I}\circ f$ with $\mathcal{I}$ as before, 
and find that $f^\perp=\frac{1}{2}(1-|f|^2)H$. 
Now let $e$  be a vector in the Lorentzian space 
$\mathbb{L}^{n+3}=\R^{n+2}\oplus\R e$ such that $\<e,e\> =-1$. Then, 
the ``inversion'' $\mathcal{T}$ with respect to the hyperbolic space
$$
\Hy_1^{n+2}(q_0)=\{p\in\mathbb{L}^{n+3}:\< p-q_0,p-q_0\>=-1\}         
$$
is given by
$$
\mathcal{T}(p)=q_0-\frac{1}{\<p-q_0,p-q_0\>}( p-q_0)
$$
with $q_0=(0,1)=e$.
If $\hat{f}=\mathcal{T}\circ f$, there is  a vector 
bundle isometry $\mathcal{\hat P}$ between the normal bundles $N_fM$ and 
$N_{\hat{f}}M$ given by
$$
\mathcal{\hat P}\mu=\mu-2\frac{\<f-q_0,\mu\>}{\<f-q_0, f-q_0\>}(f-q_0)
$$
such that shape operators of $f$ and $\tilde{f}$ are related by
$$
\hat{A}_{\mathcal{\hat P}\mu}=-\<f-q_0, f-q_0\> A_\mu-2\<f-q_0,\mu\>I 
$$
and the mean curvature vector field  of 
$\hat{f}=\mathcal{T}\circ f$ is given by
\bea
H_{\hat{f}}\!\!\!&=&\!\!\!-\mathcal{\hat P}(\<f-q_0,f-q_0 \> H+2(f-q_0)^\perp)\\
\!\!\!&=&\!\!\!2\mathcal{\hat P}e\\
\!\!\!&=&\!\!\!4(\hat f-e/2),
\eea
and thus $\hat{f}$ is minimal in 
$$
\Hy_{1/2}^{n+2}(e/2)=\Hy_R^{n+2}(q_0)
=\{p\in\mathbb{L}^{n+3}:\< p-q_0,p-q_0\>=-1/4\}. 
$$

Conversely, assume that $f$ is a composition of a minimal surface in a 
$\Hy_{1/2}^{n+2}(e/2)$  in the Lorentzian space 
$\mathbb{L}^{n+3}=\R^{n+2}\oplus\R e$ with a  stereographic projection onto 
$\R^{n+2}$. Then, we have 
$$
f^\perp=\frac{1}{2}(1-|f|^2)H.
$$
Consider the surface $\tilde{f}=\mathcal{I}\circ f$, where  $\mathcal{I}$ 
is the inversion with respect to a sphere with radius $R=1$ centered at 
the origin. We claim that $\tilde f$ is the dual to $f$. Its mean curvature 
is given by
$$
\tilde H= \mathcal{P}  ( |f|^2 H+2f^\perp)= \mathcal{P}  (H),
$$
where $\mathcal{P}$ is the corresponding bundle isometry from $N_fM$ 
to $N_{\tilde f}M$. Then $|\tilde H|=|H|$~and 
\bea
\tilde{f}+r^2\tilde H\!\!\!&=&\!\!\!\frac{f}{|f|^2}+r^2\mathcal{P}H\\
\!\!\!&=&\!\!\! \frac{f}{|f|^2}+ r^2 H-\frac{2r^2\<f^\perp,H\>}{|f|^2}f\\
\!\!\!&=&\!\!\!\frac{f}{|f|^2}+ r^2 H-\frac{1-|f|^2}{|f|^2}f\\
\!\!\!&=&\!\!\! f+ r^2 H.
\eea
The proof that (\ref{eqspe}) holds is the same as before, and thus 
$\tilde f$ is the dual to $f$.\qed

\section{The codimension two case}

In this section, we give a proof of the main result for codimension two 
in \cite{dt} making use of the computations of this paper.
\vspace{1,5ex}

For a surface $g\colon M^2\to\R^{4}$ consider the two possible
complex structures $\hat{J}_{\pm}$ on $N_gM$ and denote by
$\mathcal{J}_{\pm}$
the complex structure on the induced
bundle $g^*(T\R^4)$ given by
$$
\mathcal{J}_{\pm} \circ g_*=g_*\circ J\;\;\text{and}
\;\;\mathcal{J}_{\pm}|_{N_gM}=\hat{J}_{\pm}.
$$

\begin{theorem}\po\label{cod2}
Let $f\colon M^2\to\R^{4}$ be a superconformal surface free of minimal
and umbilical points.  Then its surface of centers $g$  is minimal and
\be\label{para}
f=g+\mathcal{J}_{\pm}h,
\ee
where $h$ is the conjugate minimal surface to $g$.
Conversely, given a simply connected minimal surface $g$ with conjugate surface
$h$, then (\ref{para}) parametrizes a superconformal surface.
\end{theorem}

\proof First assume that $f\colon M^2\to\R^4$ is a superconformal surface
and define
$$
h= r^2(\la_2\xi_1-\la_1\xi_2).
$$
Arguing as in Lemma \ref{17}, we obtain
$$
g_{*}X_1=-\dfrac{\mu}{H}f_*(X_1)
+\dfrac{\mu}{H^2}(a_1\xi_1-a_2\xi_2),\;\;
g_{*}X_2=\dfrac{\mu}{H}f_*(X_2)-
\dfrac{\mu}{H^2}(a_2\xi_1+a_1\xi_2).
$$
Moreover,
$$
h_*X_1= g_{*}X_2\;\;\mbox{and}\;\;h_*X_2= -g_{*}X_1.
$$
This shows that   $ h_*= g_*\circ J$. In particular, $g$ is minimal and 
$h$ is the conjugate surface.
We decompose $h$ into its tangent and normal components
\be\label{deco}
h=g_*(Y)+\eta.
\ee
 From Lemma \ref{ht} we know that $Y=J\grad_g\varrho$.
Using (\ref{19}), we easily obtain
$$
\< r^2 H+g_*(JY),g_*X\>=0
$$
for any $X\in TM$. Since $\eta$ is perpendicular  to $H$, we have
$$
\<r^2H+g_*(JY),\eta\>=0.
$$
Hence, $r^2 H+g_*(JY)\in N_{g}M$  is perpendicular to $\eta$.
Moreover, it is easy to see that
$$
|r^2 H+g_*(JY) |=|\eta|.
$$
This means that
$$
r^2H+g_*(JY)=-\hat{J}_{\pm}\eta.
$$
Using (\ref{19}), (\ref{20}) and (\ref{21}), we see that the tangent
component of $r^2H$ is given by
\bea
(r^2H)^{{ g_*}(T M)}
\!\!\!&=&\!\!\!\frac{r^2}{u}g_*(\<H,g_* X_1\>X_1+\<H,g_*X_2\> X_2)\\
\!\!\!&=&\!\!\!\dfrac{\mu r^4}{u}g_* (\lambda_1a_1X_1-\lambda_1a_2X_2)\\
\!\!\!&=&\!\!\!-\dfrac{\mu r^4}{u}g_*\Big(\dfrac{X_1(H^2)}{2\mu}X_1
+\dfrac{X_2(H^2)}{2\mu}X_2\Big)\\
\!\!\!&=&\!\!\!-\dfrac{r^4}{2u}g_* (\grad H^2)=\dfrac{1}{2u}g_* (\grad r^2)\\
\!\!\!&=&\!\!\!\dfrac{1}{2} g_*  (\grad_g  r^2)\\
\!\!\!&=&\!\!\!- g_*( J Y).
\eea
Hence, we obtain 
$$
f=g+g_*(JY)+ \hat{J}_{\pm}\eta.
$$

For the converse, let $g$ be a simply connected
minimal surface with conjugate $h$.  We claim that
$f=g+\mathcal{J}_{\pm}h$ is superconformal.
We decompose $h$ as in (\ref{deco}).
Differentiating with respect to $X\in TM$ and using that $h_*=g_*\circ J$ 
yields
$$
JX=\n_XY-A_\eta X,\;\; \alpha_g(X,Y)=-\nap_X \eta.
$$
Then, from $f= g+g_*(JY)+ \hat{J}_{\pm}\eta$ we find
$$
f_*X= g_*(J\circ A_\eta X- A_{\hat{J}_{\pm}\eta}X)+\alpha_g (X,JY)-
\hat{J}_{\pm}\alpha_g (X,Y).
$$
Moreover, for any $X_1, X_2 \in TM$, we have
\bea
\<f_*X_1,f_*X_2\>\!\!\!&=&\!\!\!-(\det A_\eta 
+\det A_{\hat{J}_{\pm}\eta})\<X_1,X_2\>_g
\!-\< J A_\eta X_1,A_{\hat{J}_{\pm}\eta }X_2\>\!-\< A_{\hat{J}_{\pm}\eta
}X_1,JA_\eta X_2 \>\\
\!\!\!&+&\!\!\!\<\alpha_g(X_1,Y),\alpha_g(X_2,Y)\>
+\<\alpha_g(X_1,JY),\alpha_g(X_2,JY)\>\\
\!\!\!&-&\!\!\!\<\alpha_g(X_1,JY),\hat{J}_{\pm}\alpha_g(X_2,Y)\>
-\<\alpha_g(X_2,JY),\hat{J}_{\pm}\alpha_g (X_1,Y) \>.
\eea
The Gaussian curvature $K$ and the normal curvature $K^\perp$ of $g$ satisfy
$$
\<JA_\eta X_1,A_{\hat{J}_{\pm}\eta}X_2\>+\<JA_\eta X_2,A_{\hat{J}_{\pm}\eta
}X_1\>
=|\eta|^2K^\perp \<X_1,X_2\>_g,
$$
$$
\<\alpha_g(X_1,JY),\hat{J}_{\pm}\alpha_g(X_2,Y)\>
+\<\alpha_g(X_2,JY),\hat{J}_{\pm}\alpha_g(X_1,Y)\>
=|Y|^2K^\perp\<X_1,X_2\>_g,
$$
and
$$
\<\alpha_g (X_1,Y),\alpha_g (X_2,Y)\> +\<\alpha_g(X_1,JY),\alpha_g(X_2,JY)\>
=-|Y|^2K \<X_1, X_2\>_g.
$$
Thus, the induced metric of $f$ is given by
$$
\<\,,\,\>_f=-|h|^2(K+K^\perp)\<\;,\;\>_g.
$$
The normal bundle of $f$ is spanned by $h$ and $\mathcal{J}_{\pm}h.$ Moreover,
\bea
\nab_Xf_*V\!\!\!&=&\!\!\!g_*\big(\n_X JA_\eta V-\n_X A_{\hat{J}_{\pm}\eta}V
-A_{\alpha_g(JY,V)}X+A_{\hat{J}_{\pm}\alpha_g (Y,V)}X\big)\\
\!\!\!&+&\!\!\!\alpha_g(JX,A_\eta V)-\alpha_g(X,A_{\hat{J}_{\pm}\eta}V)
+\nap_X \alpha_g (JY,V)-\hat{J}_{\pm}\nap_X\alpha_g(Y,V).
\eea
Since $g$ is minimal, we have $J\circ A_\nu=-A_\nu\circ J$ for any $\nu\in N_gM$
and $A_\nu\partial_z\in T''M$ for any  local complex coordinate $z$.
Since the (1,1)-part of the second fundamental form of $g$ vanishes, we obtain
$$
\alpha_g (\partial_z,A_\nu \partial_z)=0.
$$
 We now easily see that
\bea
\nab_{\partial_z}f_*\partial_z
\!\!\!&=&\!\!\!g_*\big(-i\n_{\partial_z}A_\eta \partial_z
-\n_{\partial_z} A_{\hat{J}_{\pm}\eta} \partial_z
-iA_{\alpha_g(Y,\partial_z)}\partial_z
+A_{\hat{J}_{\pm}\alpha_g(Y,\partial_z)}\partial_z\big)\\
\!\!\!&+&\!\!\!i\nap_{\partial_z} \alpha_g(Y,\partial_z)
-\hat{J}_{\pm}\nap_{\partial_z}\alpha_g(Y,\partial_z).
\eea
Thus, the second fundamental form of $f$ satisfies
\bea
\< \alpha_f (\partial_z,  \partial_z ), h \>
\!\!\!&=&\!\!\! \<\nab_{\partial_z} f_*{\partial_z} , h \>\\
\!\!\!&=&\!\!\!-i\<\n_{\partial_z} A_\eta {\partial_z},Y\>
-\<\n_{\partial_z}A_{\hat{J}_{\pm}\eta} {\partial_z},Y\>
-i\<a_{\alpha_g (\partial_z, Y)}{\partial_z},Y\>\\
\!\!\!&+&\!\!\!\< A_{\hat{J}_{\pm}\alpha_g(\partial_z,Y)}{\partial_z},Y\>
+i\<\nap_{\partial_z}\alpha_g(\partial_z,Y),\eta\>
-\<\hat{J}_{\pm}\nap_{\partial_z}\alpha_g(\partial_z,Y),\eta\>
\eea
and
\bea
\<\alpha_f(\partial_z,\partial_z),\mathcal{J}_{\pm}h\> \!\!\!&=&\!\!\!
\<\n_{\partial_z}A_\eta{\partial_z},Y\>
-i\<\n_{\partial_z}A_{\hat{J}_{\pm}\eta}{\partial_z},Y\>
+\<a_{\alpha_g(\partial_z,Y)}{\partial_z},Y\>\\
\!\!\!&+&\!\!\! i\<a_{\hat{J}_{\pm}\alpha_g(\partial_z,Y)}{\partial_z},Y\>
+\< \nap_{\partial_z}\alpha_g(\partial_z,Y),\hat{J}_{\pm}\eta\>
-\<\nap_{\partial_z}\alpha_g(\partial_z,Y),\eta\>.
\eea
Hence,
$$
\<\alpha_f(\partial_z,\partial_z),\mathcal{J}_{\pm}h\>
=i\<\alpha_f(\partial_z,\partial_z),h\>.
$$
This means that $\alpha_f(\partial_z,\partial_z)$ isotropic, and thus $f$ is
superconformal. \qed

\vspace{.5in} {\renewcommand{\baselinestretch}{1}
\hspace*{-20ex}\begin{tabbing} \indent\= IMPA -- Estrada Dona Castorina, 110
\indent\indent\= Univ. of Ioannina -- Math. Dept. \\
\> 22460-320 -- Rio de Janeiro -- Brazil  \>
45110 Ioannina -- Greece \\
\> E-mail: marcos@impa.br \> E-mail: tvlachos@uoi.gr
\end{tabbing}}

\begin{thebibliography}{lbl}

\bibitem{bl} W. Blaschke, Vorlesungen \"uber Differentialgeometrie und geometrische Grundlagen 
von Einsteins Relativit\"atstheorie, III: Differentialgeometrie der Kreise und Kugeln, Springer,
Berlin, 1929.

\bibitem{BPW} J. Bolton, F. Pedit and L. Woodward, \emph{Minimal surfaces and
the affine Toda field model}. J. Reine Angew. Math. \textbf{459} (1995), 119--150.

\bibitem{Ch} C. C. Chen, \emph{The generalized curvature
ellipses and minimal surfaces}. Bull. Acad. Sinica \textbf{11} (1983), 329--336.

\bibitem{DG2} M. Dajczer and D. Gromoll, \emph{The Weierstrass 
representation for complete minimal real Kaehler
submanifolds}. Invent. Math. \textbf{119} (1995), 235--242.

\bibitem{br} R. Bryant, \emph{ Conformal and minimal immersions of compact surfaces
into the $4$-sphere.} J. Diff. Geom. {\bf 17} (1982), 455--473.

\bibitem{dt0} M. Dajczer and R. Tojeiro, \emph{Commuting Codazzi tensors and the Ribaucour 
Transformation for submanifolds.} Result. Math. {\bf 44} (2003), 258--278.

\bibitem{dt} M. Dajczer and R. Tojeiro, \emph{All superconformal surfaces in $R^4$ in 
terms of minimal surfaces.} Math. Z. {\bf 261} (2009), 869--890.

\bibitem{dw} J. Dorfmeister and P. Wang, \emph{Willmore surfaces in $\Sf^{n+2}$ by 
the loop group method: generic cases and some examples.} http://arxiv.org/pdf/1305.2514.pdf

\bibitem{ej} N. Ejiri, \emph{Willmore surfaces with a duality in $\Sf^n(1)$.} 
Proc. London Math. Soc. {\bf 57} (1988), 383-416.

\bibitem{gr} I. Guadalupe and L. Rodr\'{\i}guez, \emph{Normal curvature of surfaces 
in space forms.}
Pacific J. Math. {\bf 106}  (1983), 95--103.

\bibitem{hj} U. Hertrich-Jeromin,  \emph{Introduction to
M\"obius differential geometry}, London Math. Lect. Notes Series, vol. 300, 
Cambridge Univ. Press, Cambridge, 2003.

\bibitem{ma}  X. Ma, \emph{Isothermic and $S$-Willmore surfaces as solutions to a
problem of Blaschke.}  Results Math.  {\bf 48}  (2005),  301--309.

\bibitem{Mi} R. Miyaoka, \textit{The family of isometric superconformal
harmonic maps and the affine Toda equations}, J. Reine Angew. Math. 
\textbf{481} (1996), 1-25.

\bibitem{mw1} C. Moore and E. Wilson, \emph{A general theory of surfaces.} 
J. Nat. Acad Proc.  {\bf 2} (1916), 273--278.

\bibitem{mw2} C. Moore and E. Wilson, \emph{Differential geometry of
two-dimensional surfaces in hyperspaces.} Proc. of the Academy of
Arts and Sciences, {\bf 52} (1916), 267--368.

\bibitem{ro} B. Rouxel, \emph{Harmonic spheres of a submanifold in Euclidean space.}
Proc. of the 3rd Congress of Geometry, Thessaloniki (1991), 357--364.

\bibitem{th} G. Thomsen, \emph{\"Uber konforme Geometrie I: Grundlagen der 
konformen Fl\"achentheorie.} Hamb. Math. Abh. {\bf 3} (1923), 31--56.

\bibitem{Vl} Th. Vlachos, \emph{Minimal surfaces, Hopf differentials and the
Ricci condition}. Manuscripta Math. \textbf{126} (2008), 201--230.
\end{thebibliography}
\end{document}